\documentclass[12pt]{amsart}
\usepackage{amssymb}
\usepackage{geometry}

\theoremstyle{definition}

\theoremstyle{remark}

\numberwithin{equation}{section}

\input{tcilatex}

\begin{document}
\title{K-Ring of the Classifying Space of the Symmetric Group on Four Letters%
}
\author{Mehmet K\i rdar}
\address{Department of Mathematics, Faculty of Arts and Science, Nam\i k
Kemal University, Tekirda\u{g}, Turkey}
\email{mkirdar@nku.edu.tr}
\date{28 October, 2013}

\begin{abstract}
We describe $K(BS_{4})$ and make a connection of the order of the bundle
induced from the standart representation over the four dimensional skeleton
of $BS_{4}$ with the stable homotopy group $\pi _{3}^{s}=Z_{24}$ explaining
the reasons of this connection by pulling this bundle over lens spaces.
\end{abstract}

\maketitle







\section{Introduction}

$K$-rings of the classifying spaces of cyclic groups, so-called lens spaces,
are well-known, \cite{Kirdar-Hacettepe}. In \cite{Kirdar} and \cite%
{Kirdar-Ozdemir}, we described the $K$-ring of the classifying space of the
dihedral and the generalized quaternion group in a naive way.

A natural continuation of the problem seems to be to study the $K$-rings of
the classifying spaces of symmetric groups. Since, the representation theory
of symmetric groups is very complicated, this problem turns out to be
interesting even for the first non-trivial case which is the classifying
space of the symmetric group on 4 letters. In this note, we will describe $%
K(BS_{4})$ by making connections with the integral cohomology of $S_{4},$
which is described in \cite{Thomas} via the Atiyah-Hirzebruch Spectral
Sequence (AHSS) of the $K$-ring.

We will also make a fancy connection of the $\widetilde{K}$-order of the
bundle induced from the standard representation over the four dimensional
skeleton of $BS_{4}$ with the stable homotopy group $\pi _{3}^{s}=Z_{24}.$

\section{Representations}

There are 5 conjugacy classes of $S_{4}$: $C_{1}=\left\{ 1\right\} ,$ $C_{2}=%
\overline{(12)},$ $C_{3}=\overline{(123)},$ $C_{4}=\overline{(1234)},$ $%
C_{5}=\overline{(12)(34)}.$ And thus, there are 5 irredicuble
representations: The trivial one dimensional representation $1,$ the sign
representation $d_{1}$ which is also one dimensional, the two dimensional
representation $d_{2},$ the three dimensional representation $d_{3}$ which
takes value $+1$ on the conjugacy class $C_{2},$ and the three dimensional
representation $d_{1}d_{3}.$ You can see the character table of $S_{4}$ in 
\cite{Thomas}.

The relations of the representation ring $R(S_{4})$ are: 
\begin{equation*}
d_{1}^{2}=1,\text{ }d_{2}^{2}=1+d_{1}+d_{2},\text{ }%
d_{3}^{2}=1+d_{2}+d_{3}+d_{1}d_{3},\text{ }d_{2}d_{3}=d_{3}+d_{1}d_{3},\text{
}d_{1}d_{2}=d_{2}.
\end{equation*}%
Hence, any representation of $S_{4}$ can be written as an integral linear
combination of the 5 irreducible representations described above.

From the above relations it is clear that the representation ring $R(S_{n})$
is generated by $d_{1}$ and $d_{3}$ since $%
d_{2}=d_{3}^{2}-1-d_{3}-d_{1}d_{3}.$ In fact, this observation is true for
any symmetric group: Hooks generate the representation ring of any symmetric
group, \cite{Marin}.

Let $res$ denote the restriction operator for a representation to a
subgroup. We have the following observations for the restriction of the
representations of $S_{4}$ which will be useful to understand the
filtrations of the corresponding vector bundles on the AHSS of the $K$-ring.

Let $Z_{2}$ be the cyclic subgroup generated by $(12),$ then $%
res(d_{1})=\eta ,$ $res(d_{2})=2$ and $res(d_{3})=1+2\eta $ where $\eta $
denotes the tautological one dimensional complex representation of $Z_{2}.$

Let $Z_{3}$ be the cyclic subgroup of $S_{4}$ generated by $(123),$ then $%
res(d_{1})=1,$ $res(d_{2})=\eta +\eta ^{2}$ and $res(d_{3})=1+\eta +\eta
^{2} $ where $\eta $ denotes the tautological one dimensional complex
representation of $Z_{3}.$

Let $Z_{4}$ be the cyclic subgroup of $S_{4}$ generated by $(1234),$ then $%
res(d_{1})=\eta ^{2},$ $res(d_{2})=1+\eta ^{2}$ and $res(d_{3})=\eta +\eta
^{2}+\eta ^{3}$ where $\eta $ denotes the tautological one dimensional
complex representation of $Z_{4}.$

\section{Cohomology}

The integral cohomology of $S_{4}$ is generated by 4 elements $%
a_{2},a_{3},a_{4},b_{4}$ with dimensions $\dim a_{2}=2,$ $\dim a_{3}=3,$ $%
\dim a_{4}=4$ and $\dim b_{4}=4$ indicated by subscripts. And the cohomology
ring is described in \cite{Thomas} in the following way:%
\begin{equation*}
H^{\ast }(S_{4})=Z[a_{2},a_{3},a_{4},b_{4}]/I
\end{equation*}%
where $I$ is the ideal generated by the elements $2a_{2},$ $2a_{3},$ $%
4a_{4}, $ $3b_{4}$ and, something weird, by the elements $%
a_{2}a_{3}^{2j}-a_{2}^{j+1}(a_{4}+a_{2}^{2})^{j}$ for all $j\geq 1.$

The relations $a_{2}a_{3}^{2j}=a_{2}^{j+1}(a_{4}+a_{2}^{2})^{j}$ for all $%
j\geq 1$ and also the whole information about the 2-primary part of the
cohomology are obtained by means of the cohomology of the dihedral subgroup $%
D_{8}$ mentioned in the previous section through the group inclusion
homomorphism $i:D_{8}\rightarrow S_{4}$ by Thomas and for details we advise
the reader to take a look at \cite{Thomas}.

Note also that 2-primary part and 3-primary part of the cohomology are
totally separated as it is immediate to deduce the relations $a_{2}b_{4}=0,$ 
$a_{3}b_{4}=0,$ $a_{4}b_{4}=0$ from the relations $2a_{2}=0,$ $2a_{3}=0,$ $%
4a_{4}=0,$ $3b_{4}=0.$

Since the $K$-ring of a topological space is only related to the even
dimensional part of the integral cohomology through the Atiyah-Hirzebruch
spectral sequence, we will extract only the even dimensional cohomology
groups and they can be tabulated like $H^{2j}(BS_{4})=$

$%
\begin{array}{l}
Z \\ 
Z_{2}(a_{2}) \\ 
Z_{2}(a_{2}^{2})\oplus Z_{4}(a_{4})\oplus Z_{3}(d_{4}) \\ 
Z_{2}(a_{2}^{3})\oplus Z_{2}(a_{2}a_{4})\oplus Z_{2}(a_{3}^{2}) \\ 
Z_{2}(a_{2}^{4})\oplus Z_{2}(a_{2}^{2}a_{4})\oplus Z_{4}(a_{4}^{2})\oplus
Z_{3}(d_{4}^{2}) \\ 
Z_{2}(a_{2}^{5})\oplus Z_{2}(a_{2}^{3}a_{4})\oplus
Z_{2}(a_{2}a_{4}^{2})\oplus Z_{2}(a_{3}^{2}a_{4}) \\ 
\oplus _{i=0}^{3k-1}Z_{2}(a_{2}^{6k-2i}a_{4}^{i})\oplus
_{i=0}^{k-1}Z_{2}(a_{3}^{4k-4i}a_{4}^{3i})\oplus Z_{4}(a_{4}^{3k})\oplus
Z_{3}(d_{4}^{3k}) \\ 
\oplus _{i=0}^{3k}Z_{2}(a_{2}^{6k+1-2i}a_{4}^{i})\oplus
_{i=0}^{k-1}Z_{2}(a_{3}^{4k-4i-2}a_{4}^{3i+2}) \\ 
\oplus _{i=0}^{3k}Z_{2}(a_{2}^{6k+2-2i}a_{4}^{i})\oplus
_{i=0}^{k-1}Z_{2}(a_{3}^{4k-4i}a_{4}^{3i+1})\oplus Z_{4}(a_{4}^{3k+1})\oplus
Z_{3}(d_{4}^{3k+1}) \\ 
\oplus _{i=0}^{3k+1}Z_{2}(a_{2}^{6k+3-2i}a_{4}^{i})\oplus
_{i=0}^{k}Z_{2}(a_{3}^{4k+2-4i}a_{4}^{3i}) \\ 
\oplus _{i=0}^{3k+1}Z_{2}(a_{2}^{6k+4-2i}a_{4}^{i})\oplus
_{i=0}^{k-1}Z_{2}(a_{3}^{4k-4i}a_{4}^{3i+2})\oplus Z_{4}(a_{4}^{3k+2})\oplus
Z_{3}(d_{4}^{3k+2}) \\ 
\oplus _{i=0}^{3k+2}Z_{2}(a_{2}^{6k+5-2i}a_{4}^{i})\oplus
_{i=0}^{k}Z_{2}(a_{3}^{4k+2-4i}a_{4}^{3i+1})%
\end{array}%
\begin{array}{l}
\text{for }j=0 \\ 
\text{for }j=1 \\ 
\text{for }j=2 \\ 
\text{for }j=3 \\ 
\text{for }j=4 \\ 
\text{for }j=5 \\ 
\text{for }j=6k \\ 
\text{for }j=6k+1 \\ 
\text{for }j=6k+2 \\ 
\text{for }j=6k+3 \\ 
\text{for }j=6k+4 \\ 
\text{for }j=6k+5%
\end{array}%
$ where $k\geq 1$ and $Z_{n}(g)$ denotes the cyclic group of order $n$ with
the generator $g$.

Thomas, in \cite{Thomas}, indicates that the even cohomology ring $%
H^{even}(S_{n})$ is generated by Chern characters of the standard
representation for any $n.$ We now expect that the standard representation
will play a central role in $K$-ring as well.

\section{K-Ring}

Let $d_{1},$ $d_{2},$ $d_{3}$ denote the induced vector bundles over the
classifying space $BS_{4}.$ Since, the corresponding representations
generate the representation ring, these vector bundles generate $K(BS_{4})$
due to Atiyah-Segal Completion Theorem. We set reduced elements $v=d_{1}-1,$ 
$\delta =d_{2}-2,$ $\phi =d_{3}-3$ in $\widetilde{K}(BS_{4}).$ Now, the
relations in the representation ring transform to the set of relations:%
\begin{eqnarray*}
2v &=&-v^{2} \\
3\delta +\delta ^{2} &=&v \\
4\phi +\phi ^{2} &=&3v+\delta +v\phi \\
\delta \phi &=&3v+v\phi -3\delta \\
v\delta &=&v^{2}
\end{eqnarray*}

There is an AHSS which looks like%
\begin{equation*}
E_{2}^{p,-p}=H^{p}(S_{4};\widetilde{K}(S^{p}))\Rightarrow K(BS_{4})
\end{equation*}%
on the main diagonal of the second page and which converges to $K(BS_{4}).$
The groups $E_{\infty }^{p,-p}$ on the limit page are called filtrations and
they are generated by the reduced vector bundles which are trivial on $(p-1)$%
-th skeleton of $BS_{4},$ but non-trivial on the $p$-th skeleton. In
particular, the generators of the filtrations are some particular elements
in $\widetilde{K}(BS_{4}).$ Note that these elements are related to the
cohomology of $BS_{4}$ since they are limits of the elements on the second
page. These classes probably corresponds to Chern characters of these
bundles, but we don't go on this direction further in this note.

Since $\widetilde{K}(S^{p})=Z$ when $p$ is even and zero otherwise, we have
only even dimensional filtrations $E_{\infty }^{2j,-2j}$ on the second page
of the AHSS. Because of that, only even dimensional cohomology ring $%
H^{even}(S_{4})$ plays a role for $K(BS_{4}).$

Firstly, let us explain the filtration $E_{\infty }^{2,-2}.$ We claim that
AHSS collapses on second page at $(2,-2)$ and $E_{\infty }^{2,-2}$ is $Z_{2}$
and is generated by $v$. Let us pick a subgroup $Z_{2}$ in $S_{4}$, say
generated by $(12).$ Then the subgroup inclusion $Z_{2}\rightarrow S_{4}$
induces ring homomorphism $\widetilde{K}(BS_{4})\rightarrow \widetilde{K}%
(BZ_{2})$ and also group homomorphisms $E_{\infty
}^{p,-p}(BS_{4})\rightarrow $ $E_{\infty }^{p,-p}(BZ_{2})$ on AHSSs by
naturality. Under the former map, since it corresponds to restriction of the
representation, $v$ maps to $v$ where the last same letter denotes the
tautological generator of $\widetilde{K}(BZ_{2}).$ We recall that $v\in 
\widetilde{K}(BZ_{2})$ generates $E_{\infty }^{2,-2}(BZ_{2})$ of the AHSS of 
$\widetilde{K}(BZ_{2}).$ Therefore, $v$ generates $E_{\infty
}^{2,-2}(BS_{4}).$ We also note that $v$ corresponds to the first Chern
character of the standard representation although it is the reduction of the
sign representation.

The same kind of arguments show that $v^{j}=(-2)^{j-1}v$ generates a $Z_{2}$
summand on the filtration $E_{\infty }^{2j,-2j}$ for all $j\geq 1$ and these
summands correspond to the summands $Z_{2}(a_{2}^{j})$ on $E_{2}^{2j,-2j}.$

Next, we will deal with $E_{\infty }^{4,-4}$. As we said the $Z_{2}$ summand
of $E_{\infty }^{4,-4}$ is generated by $v^{2}.$ Now, we should explain its
remaining $Z_{12}$ summand for which we claim that it survives.

Let us pick a subgroup $Z_{4}$ in $S_{4}$, say generated by $(1234).$ Then
the subgroup inclusion $Z_{4}\rightarrow S_{4}$ induces ring homomorphism $%
\widetilde{K}(BS_{4})\rightarrow \widetilde{K}(BZ_{4})$ and also group
homomorphisms $E_{\infty }^{p,-p}(BS_{4})\rightarrow $ $E_{\infty
}^{p,-p}(BZ_{4})$ on AHSSs by naturality. At this point we observe that the
pull-back of $\phi $ to $BZ_{4}$ lives on $E_{\infty }^{2,-2}(BZ_{4})$ so
that it can not be a generator of $E_{\infty }^{4,-4}$ of the AHSS of our
ring.

Because of that, we set $x=\phi +v$ and now $x$ should be a generator of $%
E_{\infty }^{4,-4},$ since it creates no conflict with the restrictions to
cyclic subgroups $Z_{4}$ and $Z_{3}$ as well: $x$ pulls back to a bundle in
the form $\mu ^{2}+...$ over both $BZ_{4}$ and $BZ_{3}$ where $\mu =\eta -1$
is the reduction of the tautological line bundle. Thus it must have order $%
12 $ in $E_{\infty }^{4,-4}.$ Therefore, the AHSS collapses on second page
at the coordinate $(4,-4)$ too. But, it will not be the case on the
dimension 6.

Really, by the similar elementary explanations, the higher filtrations are
observed as below:

$E_{\infty }^{2,-2}=Z_{2}(v)$

$E_{\infty }^{4,-4}=Z_{2}(v^{2})\oplus Z_{12}(x)$

$E_{\infty }^{6,-6}=Z_{2}(v^{3})\oplus Z_{2}(v\phi )$ (AHSS is not
collapsing from hereafter!)

$E_{\infty }^{8,-8}=Z_{2}(v^{4})\oplus Z_{2}(v^{2}\phi )\oplus Z_{12}(x^{2})$

$E_{\infty }^{10,-10}=Z_{2}(v^{5})\oplus Z_{2}(v^{3}\phi )\oplus Z_{2}(v\phi
^{2})$

$E_{\infty }^{12,-12}=Z_{2}(v^{6})\oplus Z_{2}(v^{4}\phi )\oplus
Z_{2}(v^{2}\phi ^{2})\oplus Z_{12}(x^{3})$

.

.

.

etc.

We observe that only the parts of the even cohomology related to the class $%
a_{3}\in H^{3}(S_{4})$ are not surviving on the AHSS; they vanish. The other
parts of the even cohomology stay as they are. Hence, we have

\textbf{Theorem 1:\ }%
\begin{equation*}
K(BS_{4})=Z[v,\phi ]\diagup \left( 
\begin{array}{c}
2v+v^{2} \\ 
12\phi +7\phi ^{2}+\phi ^{3}-4v-v\phi \\ 
24\phi +26\phi ^{2}+9\phi ^{3}+\phi ^{4} \\ 
2v\phi -8v-24\phi ^{2}-14\phi ^{3}-2\phi ^{4}+v\phi ^{2}%
\end{array}%
\right)
\end{equation*}

As we expected, the $K$-ring is just a change of variables in the $R$-ring.

\section{A\ Little Homotopy}

Third stable homotopy group of spheres is $\pi _{3}^{s}=Z_{24}$ and it is
the image of the $J$-homomorphism 
\begin{equation*}
J:\widetilde{KO}(S^{4})\rightarrow \pi _{3}^{s}.
\end{equation*}%
\medskip\ Its generator is 
\begin{equation*}
S^{4}\overset{1}{\rightarrow }BO\overset{Bi}{\rightarrow }B(Z\times
BS_{\infty }^{+}).
\end{equation*}

Now, it follows from the previous section, that the order of the map%
\begin{equation*}
\phi :BS_{4}^{(4)}\rightarrow BU
\end{equation*}%
is $24$ where $BS_{4}^{(4)}$ is the four dimensional skeleton of $BS_{4}$.
We believe that it is not a coincidence that the order of $\phi $ in $%
\widetilde{K}(BS_{4}^{(4)})$ is the same as the order of $\pi _{3}^{s}$.

We can explain the reasons a little bit behind our claim by pulling the map $%
\phi $ over some lens spaces which is our fundamental trick. As we said
before, we have natural subgroups which induce projections $%
BZ_{3}\rightarrow BS_{4}$ and $BZ_{4}\rightarrow BS_{4}$ and the pull-backs
of $\phi $ over these lens spaces by these projections are equal to $\mu ^{2}
$ and $2\mu +h.o.t.$ respectively where $\mu $ is the reduction of the
corresponding Hopf bundle. On the other hand, $r(\mu ^{2})=2w+w^{2}$ where $r
$ is the relaification map and $r(2\mu )=2w.$ Now, the $KO$-order, and the $J
$-order as well, of $2w$ over the skeletons $BZ_{3}^{(4)}$ and $BZ_{4}^{(4)}$
are $3$ and $4$ respectively. These should correspond to the local $J$-order
of the generator of the group $\widetilde{KO}(S^{4})$ at prime $3$ and to
the half of the local $J$-order of that generator at prime $2,$
respectively. We guess, we should write down $KO$-ring of $BS_{4}$ to
understand that $2$-local behaviour since realification seems to change the
order at prime $2$.

Having this connection, we will make the following expected conjectures.

For all $n$ and in the limit $n\rightarrow \infty ,$ let $\mathbf{\phi }$
denote the family of the standard generators of $\widetilde{K}(BS_{n})$
which are induced from the standard representations. We claim that $\mathbf{%
\phi }$ detect the $\func{Im}J$ part of $\pi _{\ast }^{s}$. This conjecture
seems to be quite possibly true due to Hegenbarth's work, \cite{Hegenbarth},
on the $k$-theory of the classifying space of the infinite symmetric group
where his $\func{mod}$ $p$ generators suggest that connection.

A stronger and much more fancy conjecture would be that

\textbf{Conjecture 2: }$\mathbf{\phi }$ detect $\pi _{\ast }^{s}$.

In other words, we claim that representations of symmetric groups are saying
quite a lot about the homotopy groups of spheres. We note that it becomes
incredibly complicated even to write the polynomial relations in $K(BS_{n})$
when $n\rightarrow \infty ,$ yet alone the unimaginable geometric work
behind. $n=5$ stays as a messy exercise.

\end{document}